\documentclass[twocolumn]{autart}    
\usepackage{graphicx}
\usepackage{xcolor}
\usepackage{epstopdf}
\usepackage{amsmath,amssymb,amsfonts,bm}
\usepackage{natbib}
\usepackage{subfigure}

\usepackage{algorithm}
\usepackage{algpseudocode}

\allowdisplaybreaks
\def \qed {\hfill \vrule height6pt width 6pt depth 0pt}

\begin{document}
	\begin{frontmatter} 
	\title{Minimax Q-learning Control for Linear Systems Using the Wasserstein Metric} 
	\thanks[footnoteinfo] {This research was supported by National Natural Science Foundation of China under Grant no. 62033006, and Tsinghua-Foshan Innovation Special Fund (TFISF).} 
	\author{Feiran~Zhao},
	\ead{zhaofr18@mails.tsinghua.edu.cn}
	\author{Keyou~You\corauthref{cor}}
	\corauth[cor]{Corresponding author} 
	\ead{youky@tsinghua.edu.cn}
	\address{Department of Automation, Beijing National Research Center for Information Science and Technology, Tsinghua University, Beijing 100084, China}
\begin{abstract}
Stochastic optimal control usually requires an explicit dynamical model with probability distributions, which are difficult to obtain in practice. In this work, we consider the linear quadratic regulator (LQR) problem of unknown linear systems and  adopt a Wasserstein penalty to address the distribution uncertainty of additive stochastic disturbances. By constructing an equivalent deterministic game of the penalized LQR problem, we propose a Q-learning method with convergence guarantees to learn an optimal minimax controller. 
\end{abstract}
\end{frontmatter}

\section{Introduction}

Stochastic control is a well-studied framework for dynamical systems with known probability distributions of uncertain variables \citep{aastrom2012introduction}, which is usually difficult to obtain in practice. Directly using an approximated distribution is not always reliable and may even lead to system instability~\citep{nilim2005robust}. 
To tackle distribution errors, it is natural to use the idea of distributionally robust (DR) optimization with an appropriate distribution metric~\citep{gao2016distributionally,wiesemann2014distributionally}. Since the Wasserstein metric can well combine the distribution prior knowledge and sampled data~\citep{chen2021distributionally}, it has been recently  adopted to the stochastic control. Unfortunately, they require an explicit dynamical model of the control system~\citep{yang2021w, kim2020minimax, kim2021distributional}. 

Reinforcement learning (RL), as a learning-based method to approximately solve Markov decision processes (MDPs), has achieved tremendous progresses in the continuous control~\citep{lillicrap2016continuous}. Since the RL based controller is usually trained over a simulator, it may fail to generalize to practical applications if the simulator cannot generate samples with probability distributions that correctly reflect the physical system. To tackle it, robust RL algorithms are devised to enhance the robustness of the resulting control policy~\citep{pinto2017robust,abdullah2019wasserstein}. For example, \citet{abdullah2019wasserstein} uses the idea of over-fitting to training environments by minimizing a long-term cost function over a Wasserstein ambiguity set. However, they only focus on the generic MDPs and may require a large number of samples, and the stability of the closed-loop system cannot be guaranteed. 

This work takes advantages of both the DR optimization and RL to address the distribution uncertainty in an empirical measure of disturbance samples for unknown linear stochastic systems. Without explicitly identifying the explicit linear model, we aim to solve the minimax linear quadratic regulator (LQR) problem with a Wasserstein penalty via Q-learning. The major challenge lies in that the Q-function cannot be directly evaluated as it involves an expectation over the random disturbance, which is resolved by proposing an equivalent  {\em deterministic} zero-sum game. Then, {we apply \citet[Section 3]{al2007model} to design an efficient Q-learning method with convergence guarantees to learn an optimal minimax controller. }

The rest of the presentation is structured as follows. In Section \ref{sec:background}, we describe the minimax control problem with a Wasserstein penalty and provide its closed-form solution. Section \ref{sec:algorithm} designs a Q-learning algorithm to learn the minimax controller by establishing an equivalent deterministic game. Section \ref{sec:experiment} validates the convergence of the Q-learning algorithm via simulations. We draw some concluding remarks in Section \ref{sec:conclusion}.

\section{Wasserstein minimax control}\label{sec:background}
In this section, we formulate the LQR problem with a Wasserstein penalty and provide its closed-form solution.
\subsection{Problem formulation}
We consider a time-invariant linear stochastic system
\begin{equation}
\label{equ:sys}
x_{k+1} = Ax_k + Bu_k + Ew_k
\end{equation}
where $(A,B,E)$ are unknown model parameters, $x_k \in \mathbb{R}^n$ is the state and $u_k \in \mathbb{R}^m$ is the control. The disturbance $\{w_k\}\subseteq \mathbb{R}^d $ is i.i.d. with an unknown probability distribution $ \bf \nu $. Moreover, we have collected a finite number of disturbance samples $\{v^{j}\}_{j=1}^N$. Then, it is natural to use the empirical distribution, i.e., $
{\nu}=\frac{1}{N}\sum_{j=1}^{N}\delta_{v^{j}}
$
where the Dirac delta measure $\delta_{v^{j}}$ is concentrated at $v^{j}$. 


To measure the distance between two probability distributions, we adopt the following Wasserstein metric~\citep{gao2016distributionally}
$$
W(\mu_1, \mu_2) = \inf _{\kappa}\left(\int_{\mathbb{R}^d\times \mathbb{R}^d}      \left\|w_{1}-w_2\right\|^{2} \kappa \left(\mathrm{d} w_{1}, \mathrm{d} w_2\right)\right)^{{1}/{2}},
$$
where $\|\cdot\|$ denotes the standard $l_2$-norm, $\kappa \in \Gamma (\mu_1, \mu_2)$, and $\Gamma (\mu_1, \mu_2)$ denotes the set of the joint distributions with marginal distributions $\mu_1$ and $\mu_2$. 

To address the distributional uncertainty in ${\nu}$, we consider a zero-sum Markov game where the controller is the minimizing player and the disturbance is the maximizing player with distributions as inputs to the game. Particularly, we only focus on the deterministic policies. Let the control policy be $\{\pi_k\}$ that maps the history trajectory $h_k= \{x_0,u_0,\dots, x_k\}$ to the control input, i.e., $u_k = \pi_k(h_k)$. Let the adversarial policy (distribution) be $\{\gamma_k\}$ which maps $h_k$ and $u_k$ to a distribution $\mu_k$, i.e., $\mu_k = \gamma_k(h_k,u_k)$. 

\begin{figure}
	\begin{center}
		\includegraphics[width=70mm]{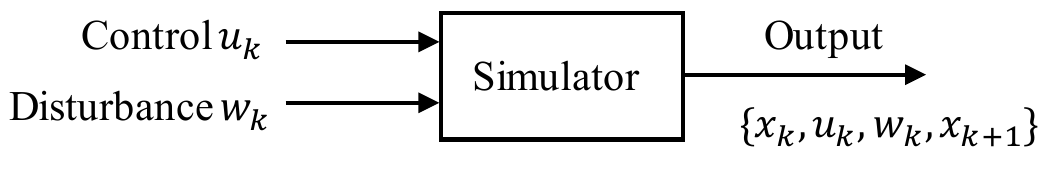}
		\caption{The simulator (a black-box to the control designer) produces the output of \eqref{equ:sys}, which together with the input serve as samples for our learning algorithm. }
		\label{pic:simu}                                
	\end{center}                                 
\end{figure}

Without identifying model $(A, B, E)$, this work aims to design a Q-learning algorithm via the simulator in Fig. \ref{pic:simu}  to solve a minimax control problem \citep{yang2021w}
\begin{equation}\label{equ:wasprob}
\min_{\{\pi_k\}} \max_{\{\gamma_k\}} \mathbb{E}^{\{\pi_k\},\{\gamma_k\}}\big[\sum_{k=0}^{\infty} \alpha^{k} c(x_k, u_k, \mu_k) \big],
\end{equation}
where $0<\alpha<1$ is a discount factor and
$c(x_k, u_k, \mu_k) = x_k^{\top}Qx_k + u_k^{\top}Ru_k - \lambda \cdot W^2\left(\mu_k, {\nu} \right)$ with $Q\ge 0$ and $R>0$. The tunable hyperparameter $\lambda>0$ reflects our confidence in the empirical distribution ${\nu}$. Specifically, the larger $\lambda$, the more we trust ${\nu}$.

\subsection{Wasserstein minimax controller}
In this work, we make the following assumptions. 
\begin{assum}\label{assume} 
	$\Lambda = BR^{-1}B^{\top} - EE^{\top}/\lambda >0$, and $(A, \sqrt{Q})$ is observable.
\end{assum}

\begin{assum}\label{assum:lambda}
	$\lambda > \underline{\lambda}$ where
	$\underline{\lambda} = \inf \{\lambda|\lambda I -\alpha E^{\top} P E > 0 \},$
	where $P$ is a positive semi-definite solution to
	\begin{equation}
	\label{def:P}
	P = H_{xx}-
	\begin{bmatrix}
	H_{xu}&H_{xw}
	\end{bmatrix}
	\begin{bmatrix}
	H_{uu} & H_{uw} \\
	* & H_{ww}
	\end{bmatrix}^{-1}	
	\begin{bmatrix}
	H_{xu}^{\top} \\
	H_{xw}^{\top}
	\end{bmatrix}
	\end{equation}
	with $H_{xx} = Q+\alpha A^{\top}PA,  H_{xu} = \alpha A^{\top}PB, H_{xw} = \alpha A^{\top}PE, H_{uu}= R+\alpha B^{\top}PB, H_{uw}= \alpha B^{\top}PE, H_{ww} = \alpha E^{\top}PE-\lambda I$.
\end{assum}

For simplicity, define $g$ and $z$ as
\begin{equation}\label{def:g}
\begin{aligned}
&\hspace{-0.285cm}g= G_x-
\begin{bmatrix}
H_{xu}&H_{xw}
\end{bmatrix}
\begin{bmatrix}
H_{uu} & H_{uw} \\
* & H_{ww}
\end{bmatrix}^{-1}
\begin{bmatrix}
G_u \\
G_w
\end{bmatrix},\\
&\hspace{-0.285cm}z =  \alpha z - \lambda \|\bar{w}\|^2 - \text{tr}\{H_{ww}^{-1} (\lambda^2\Sigma + \frac{1}{4}G_wG_w^{\top}) \}\\
&\hspace{-0.285cm}-\frac{1}{4} (G_u^{\top}-G_w^{\top}H_{ww}^{-1} H_{uw}^{\top})(H_{u u}-H_{u w}H_{w w}^{-1} H_{uw}^{\top})^{-1}\\
&\hspace{-0.285cm}\times ( G_u- H_{uw} H_{ww}^{-1}  G_w)
\end{aligned}	
\end{equation}
with $G_x = \alpha A^{\top}g, G_u = \alpha B^{\top}g, G_{w} = \alpha E^{\top}g+2\lambda \bar{w}$, $G_{w}^j = \alpha E^{\top}g+2\lambda v^{j}$. 

Via a simple extension of \citet{kim2020minimax}, we then obtain a closed-form solution to (\ref{equ:wasprob}).

\begin{lem}\label{theorem:solution}
	
	(a) The value function $V(x) = x^{\top}Px + g^{\top}x + z$ solves the following Bellman equation
	$$V(x_k) = \min_{\{\pi_k\}} \max_{\{\gamma_k\}} \mathbb{E}^{\{\pi_k\},\{\gamma_k\}}[ c(x_k, u_k, \mu_k)  + \alpha V(x_{k+1})].$$
	(b) The minimax controller for (\ref{equ:wasprob}) is stationary and has an affine structure, i.e., 
	$
	\pi_k(x) = Kx + r
	$
	where	
	\begin{equation}\label{K}
	\begin{aligned}
	&K = (H_{u u}-H_{u w} H_{w w}^{-1} H_{uw}^{\top})^{-1}(H_{u w} {H_{w w}^{-1}} H_{xw}^{\top}-H_{u x}^{\top}), \\
	&r = -\frac{1}{2}(H_{u u}-H_{u w} H_{w w}^{-1} H_{uw}^{\top})^{-1}(G_u - H_{uw} H_{w w}^{-1} G_w ).
	\end{aligned}
	\end{equation}
	(c) One of the least favorable adversarial policy to solve (\ref{equ:wasprob}) is discrete and stationary.  Specifically, let $w^j(x) = Lx + l^j$ with
	\begin{equation}\label{L}
	\begin{aligned}
	&L = (H_{w w}-H_{uw}^{\top} H_{uu}^{-1} H_{u w})^{-1}(H_{uw}^{\top} H_{uu}^{-1} H_{xu}^{\top}-H_{xw}^{\top}), \\
	&l^j = \frac{1}{2N}(H_{w w}-H_{uw}^{\top} H_{uu}^{-1} H_{u w})^{-1}( H_{uw}^{\top} H_{uu}^{-1} G_u -G^j_{w}).
	\end{aligned}
	\end{equation}
	Then, $\gamma_k(x,u)=\frac{1}{N} \sum_{j=1}^{N} \delta_{w^j}$ is a  least favorable adversarial policy. 
	\end{lem}

For a sufficiently large discount factor $\alpha$, it can be shown that the minimax controller in (\ref{K}) can also stabilize the system; see e.g., \cite{bacsar2008h}.

Computing the minimax controller in (\ref{K}) explicitly requires to know the model matrices $(A,B)$. In the sequel, we design a Q-learning algorithm to learn it via the simulator in Fig. \ref{pic:simu}.

\section{The Q-learning method and convergence}\label{sec:algorithm}

In this section, we first propose a new deterministic game which has the same minimax controller as (\ref{K}). Then, we develop a Q-learning algorithm with convergence guarantees to learn the minimax controller by solving the proposed deterministic game. 

\subsection{An equivalent deterministic game}

The Q-function $Q(x_k,u_k,\mu_k)$ of (\ref{equ:wasprob}) is given by
\begin{equation}\label{qexpectation}
\begin{aligned}
Q(x_k,u_k,\mu_k) & =  x_k^{\top}Qx_k + u_k^{\top}Ru_k - \lambda W^2(\mu_k, {\nu} ) \\
&+ \alpha \mathbb{E}_{ w_k \sim \mu_k}[
V(x_{k+1})].
\end{aligned}
\end{equation} 
Once Q-function is determined, an optimal controller can be easily obtained. However, it is ineffective to evaluate $Q(x_k,u_k,\mu_k)$ in \eqref{qexpectation} as it involves an expectation operator.  We use the Kantorovich duality \citep{gao2016distributionally} and the quadratic form of $V(x)$ to resolve it. Since 
\begin{equation}
\begin{aligned}
&\min_{u_k} \max_{\mu_k} Q(x_k,u_k,\mu_k)=\min_{u_k}  \{ x_k^{\top}Qx_k + u_k^{\top}Ru_k \\
&+  \frac{1}{N} \sum_{j=1}^{N} \max_{w}  \{ \alpha V(x_{k+1}) - \lambda  \|w-v^{j}\|^{2}  \}   \} \\
& = \min_{u_k}  \{ x_k^{\top}Qx_k + u_k^{\top}Ru_k +  \max_{w_k} \{ \alpha V(x_{k+1}) \\
&~~~- \lambda  \|w_k-\bar{w}\|^{2}  \}  \} + z_0
\end{aligned}
\end{equation}
where $z_0$ is a constant, the expectation issue in \eqref{qexpectation} can be avoided via  the following deterministic zero-sum game
\begin{equation}\label{equ:newgame}
\min_{\{\tilde{\pi}_k\}} \max_{\{\tilde{\gamma}_k\}} \sum_{k=0}^{\infty} \alpha^{k} (x_k^{\top}Qx_k + u_k^{\top}Ru_k
- \lambda \|w_k - \bar{w}  \|^2 ),
\end{equation}
where $\{\widetilde{\pi}_k\}$ is a control policy sequence in the form $u_k = \widetilde{\pi}_k(h_k)$ and $\{\widetilde{\gamma}_k\}$ is a disturbance policy sequence in the form $w_k = \widetilde{\gamma}_k(h_k,u_k)$. 

Particularly, we show in the following theorem that the minimax controller of (\ref{equ:wasprob}) and (\ref{equ:newgame}) are the same. 

\begin{thm}
	\label{theorem:certainty}
Consider the deterministic zero-sum game in (\ref{equ:newgame}).  	We obtain the following results.  

	(a) The value function $V^d(x_k) = \min_{u_k} \max_{w_k} (x_k^{\top}Qx_k + u_k^{\top}Ru_k
	- \lambda \|w_k - \bar{w}  \|^2   + \alpha V^d(x_{k+1}))$ is  given by
	\begin{equation}
	\notag
	\begin{array}{lll}
	V^d(x) & =x^{\top} P x+ g^{\top}x + z -z_0,& \\
	\end{array}
	\end{equation}
	where $P$, $g$ and $z$ are defined in (\ref{def:P}) and (\ref{def:g}).
	
	(b) The Q-function $Q(x_k,u_k,w_k) = x_k^{\top}Qx_k + u_k^{\top}Ru_k - \lambda \|w_k - \bar{w}  \|^2 + \alpha V^d(x_{k+1})$ is given by
	\begin{equation}\label{def:Q_c}
	\begin{aligned}
	Q(x,u,w) =&
	\begin{bmatrix}
	x \\
	u \\
	w
	\end{bmatrix}^{\top}
	H
	\begin{bmatrix}
	x \\
	u \\
	w
	\end{bmatrix} +
	G^{\top}
	\begin{bmatrix}
	x \\
	u \\
	w
	\end{bmatrix}
	+ s,
	\end{aligned}
	\end{equation}	
	where $s$ is a constant scalar,
	\begin{align*}
	H &=\begin{bmatrix}
	H_{xx} & H_{xu} & H_{xw} \\
	& H_{uu} & H_{uw} \\
	*&  & H_{ww}
	\end{bmatrix},
	~\text{and}~
	G^{\top} =
	\begin{bmatrix}
	G_x^{\top}, G_u^{\top}, G_w^{\top}
	\end{bmatrix}.
	\end{align*}
	(c) A minimax controller of (\ref{equ:newgame}) is given by (\ref{K}) and a maximizing policy is given by $w(x)= Lx+l,$ where $L$ is given by (\ref{L}) and
	$$l = \frac{1}{2}(H_{w w}-H_{uw}^{\top} H_{uu}^{-1} H_{u w})^{-1}(  H_{uw}^{\top} H_{uu}^{-1} G_u -G_w).
	$$ 
\end{thm}
\begin{pf}
	The proof is similar to that of \citet[Theorem 3.7]{bacsar2008h} and is omitted for saving space. \qed
\end{pf}


The major advantage of (\ref{equ:newgame}) lies in that the adversarial input is no longer a sequence of distributions.  

\subsection{Online Q-learning algorithm to solve (\ref{equ:wasprob})}
\label{sec_qlearning}
{By \eqref{qexpectation}-\eqref{def:Q_c}, we now apply the results in \citet[Section 3]{al2007model} to provide a Q-learning algorithm to learn  $Q(x,u,w)$ in (\ref{def:Q_c}) via the simulator in Fig. \ref{pic:simu}. } We first reformulate $Q(x,u,w)$ as a linear function with a parameter vector $\theta$. Let $e = [x^{\top},u^{\top} ,w^{\top}]^{\top} \in \mathbb{R}^q$ with $q = n+m+d$, $\bar{e} = [e_1^2, \cdots, e_1e_q, e_2^2, e_2e_3, \cdots, e_{q-1}e_q,e_q^2]^{\top}$, $h$ be the vector formed by stacking the columns of $H$ and then removing the redundant terms introduced by the symmetry of $H$. Then, the Q-function can be written as
\begin{equation}
\label{equ:para_q}
\begin{aligned}
Q(x,u,w |\theta) =
[h^{\top},G^{\top},s]
[\bar{e}^{\top},e^{\top},1]^{\top} = \theta^{\top} \tilde{e}
\end{aligned}
\end{equation}
with $\theta = [h^{\top},G^{\top},s]^{\top}$ and
$
\tilde{e} = [\bar{e}^{\top},e^{\top},1]^{\top}.
$

Let $d(x_k,\theta) = x_k^{\top}Qx_k + u_k^{*\top}Ru^*_k- \lambda \|w^*_k - \bar{w}\|^2 + \alpha Q(x_{k+1},u^*_{k+1},w^*_{k+1}|\theta)$,
which can be computed via the sample $\{x_k, u^*_k,w^*_k, x_{k+1}\}$ from the simulator.
Then, it follows from the Bellman equation that 
\begin{equation}\label{equ:iter}
Q(x_k,u^*_k,w^*_k|\theta)=d(x_k,\theta).
\end{equation}
Next, we find $\theta$ by designing an iterative learning algorithm to solve (\ref{equ:iter}). Suppose that at the $i$-th iteration, the parameter vector is denoted as $\theta_{i}$ and the resulting pair of policies are $u^i(x)=K_ix+r_i, w^i(x)= L_ix+l_i$. With the sampled trajectory $\{x_k^i, u^i_k, w^i_k, x_{k+1}^i\}_{k=0}^{M-1}$ of length $M$ from the simulator, $\theta_{i+1}$ is obtained by solving a least-squares problem
\begin{equation}
\label{equ:theta_update}
\begin{aligned}
	\theta_{i+1}
	&=\arg \min_\theta  \sum_{k = 0}^{M-1} |	Q(x_k^i, u^i_k, w^i_k|\theta) 
	- d(x_k^i, \theta_i)  |^2  \\
	 &=\arg \min_\theta \sum_{k = 0}^{M-1} |	{\theta}^{\top} \tilde{e}(x_k^i) - d(x_k^i, \theta_i)  |^2.
\end{aligned}
\end{equation}
Since $u^i_k$ and $w^i_k$ are linearly dependent on $[x_k^{i\top} ~~1]^{\top}$, solving (\ref{equ:theta_update}) yields an infinite number of solutions. To this end, we manually add exploration noises to the control and disturbance inputs, i.e.,
$
	{u}^i_k = K_i x_k + r_i + o_k^1,~~{w}^i_k = L_i x_k + l_i + o_k^2,
$
where $o_k^1 \sim \mathcal{N}(0, I_m)$ and $o_k^2 \sim \mathcal{N}(0, I_d)$. This ensures that if $M > \frac{1}{2}(q+1)(q+2)$,  there is a unique solution to (\ref{equ:theta_update}), i.e.,
\begin{equation}
\label{equ:ls}
\theta_{i+1} = (\sum_{k = 0}^{M-1}\tilde{e}(x_k^i)\tilde{e}(x_k^i)^{\top} )^{-1}
\sum_{k = 0}^{M-1} \tilde{e}(x_k^i) d(x_k^i,\theta_i).
\end{equation}
We provide the detailed Q-learning algorithm in Algorithm \ref{alg:q_learning} which terminates if the increment of ${\theta}_{i}$ is smaller than a user-defined threshold.
\begin{algorithm}[t]
	\caption{The Q-learning algorithm}
	\label{alg:q_learning}
	\begin{algorithmic}[1]
		\Require
		Penalty parameter $\lambda$, discount factor $\alpha$, disturbance samples $\{v^{j}\}_{j=1}^N$,
		length of a trajectory $M$.
		\Ensure
		a minimax controller $u(x)= K_{\infty}x+r_{\infty}$.
		\State  Initialize ${\theta}_0 = 0$, $K_0 = 0$, $r_0 = 0$, $L_0 = 0$ and $l_0 = 0$.
		\For{$i=0,1,\cdots $}
		\State Sample a trajectory $\{x_k^i, u^i_k, w^i_k, x_{k+1}^i\}_{k=0}^{M-1}$ \phantom{.....}with a fixed initial state $x_0$ from  the simulator.
		\State Update $\theta_{i+1}$ by (\ref{equ:ls}).
		\State Update the policy pair $(K_{i+1}, r_{i+1})$ and
		\phantom{......}$(L_{i+1}, l_{i+1})$ from $\theta_{i+1}$.
		\EndFor
	\end{algorithmic}
\end{algorithm}

For simplicity, let $W=\text{diag}(Q, R, -\lambda I)$, $c_i = [0,r_i^{\top},l_i^{\top}]^{\top}$ and $U_i = [I,K_i^{\top},L_i^{\top}]^{\top}$. Then, we prove the convergence of Algorithm \ref{alg:q_learning}.

{\begin{thm}\label{thm:convergence}
If $M > \frac{1}{2}(q+1)(q+2)$, then $\{\theta_{i}\}$ in Algorithm \ref{alg:q_learning} converges to $\theta$.
\end{thm}}
\begin{pf}
	We first show that $\theta_{i+1} = [h_{i+1}^{\top}~G_{i+1}^{\top}~s_{i+1}]^{\top}$ in (\ref{equ:ls}) is given by
	\begin{align}
	\label{equ:para_iter1}
	H_{i+1}&=
	W+\alpha \begin{bmatrix}
	A & B & E \\
	K_{i} A & K_{i} B & K_{i} E \\
	L_{i} A & L_{i} B & L_{i} E
	\end{bmatrix}^{\top}
	H_{i}\begin{bmatrix}
	A & B & E \\
	K_{i} A & K_{i} B & K_{i} E \\
	L_{i} A & L_{i} B & L_{i} E
	\end{bmatrix},\\ \label{equ:para_iter2}
	G_{i+1}^{\top}&=
	\alpha (G_i^{\top} + 2 c_i^{\top} H_i) U_i[A,B,E] + 
	[0,0,2\lambda\bar{w}^{\top}],\\\notag
	s_{i+1} &=
	\alpha (s_i + G_i^{\top} c_i + c_i^{\top} H_i c_i
	) - \lambda \|\bar{w}\|^2.
	\end{align}
	To this end, define $V_i = [\tilde{e}(x_0^i), \tilde{e}(x_1^i), \cdots, \tilde{e}(x_{M-1}^i)]$ and $Y_i = [d(x_0^i, \theta_i), d(x_1^i, \theta_i), \cdots, d(x_{M-1}^i, \theta_i)]^{\top}$. It follows from (\ref{equ:ls}) that $\theta_{i+1} = (V_iV_i^{\top})^{-1}V_iY_i$. Next, we show that $Y_i=V_i^{\top}y(\theta_i)$, where $y(\theta_{i})$ is a function of $\theta_{i}$. By the definition of $d(x_k, \theta_i)$, we obtain
	\begin{equation}\label{equ:dd}
	\begin{aligned}
	d(x_k, \theta_i) & = e(x_k)^{\top}We(x_k) + [0, 0, 2\lambda \bar{w}^{\top}]e(x_k)  - \lambda \|\bar{w}\|^2 \\
	&+ \alpha  (e(x_{k+1})^{\top}H_i e(x_{k+1}) + G_i^{\top} e(x_{k+1}) +s_i).
	\end{aligned}
	\end{equation}
	Moreover, the term $e(x_{k+1})$ in (\ref{equ:dd}) can be written as
	\begin{equation}\label{equ:e_e}
	\begin{aligned}
	e(x_{k+1})&=e\left([A,B,E]e(x_k)\right)\\
	&=e ([A,B,E] (
	U_i x_k +
	c_i
	) ) \\
	&=U_i
	[A,B,E] (
	U_i x_k +
	c_i
	) +
	c_i\\
	&=U_i
	[A,B,E]e(x_k) +
	c_i. \\
	\end{aligned}
	\end{equation}
	Inserting (\ref{equ:e_e}) into (\ref{equ:dd}), $d(x_k, \theta_i)$ can be written as a quadratic function with respect to $e(x_k)$. Hence, $d(x_k, \theta_i)$ is linear in $\tilde{e}(x_k)$ by the definition of $\tilde{e}$, i.e., $d(x_k, \theta_i) = \tilde{e}^{\top}(x_k) y(\theta_{i})$. Then, we have that $\theta_{i+1} = (V_iV_i^{\top})^{-1}V_iY_i = y(\theta_{i})$. Based on the derivation in (\ref{equ:dd}) and (\ref{equ:e_e}), it is straightforward to write the expression of $y(\theta_{i})$, which leads to (\ref{equ:para_iter1}) and (\ref{equ:para_iter2}). 
	
	Then, define $P_{i}= U_i^{\top} H_i U_i$. We note that (\ref{equ:para_iter1}) can be written as
	\begin{equation}\label{equ:18}
		\begin{aligned}
		&H_{i+1}=W+[A,B,E]^{\top}P_i[A,B,E]\\
		&\hspace{-.1cm}=\begin{bmatrix}
		\alpha A^{\top} P_i A+ Q & \alpha A^{\top} P_i B & \alpha A^{\top} P_i E \\
		&\alpha B^{\top} P_i B+ R & \alpha B^{\top} P_i E \\
		*&  & \alpha E^{\top} P_i E-\lambda I
		\end{bmatrix}.
		\end{aligned}
	\end{equation}
	Similarly, define $g_i^{\top} = (G_i^{\top} + 2 c_i^{\top} H_i)U_i$ and $z_i = s_i + G_i^{\top} c_i + c_i^{\top} H_i c_i$. Then, (\ref{equ:para_iter2}) can be written as 	
	\begin{equation}\label{equ:19}
	\begin{aligned}
	G_{i+1}^{\top} & = \left[\alpha g_i^{\top}A , \alpha g_i^{\top}B , \alpha g_i^{\top}E+2\lambda \bar{w}^{\top}  \right],\\
	s_{i+1}& = \alpha z_i - \lambda \|\bar{w}\|^2.
	\end{aligned}
	\end{equation}

	By the definition of $P_i$, we have that
	\begin{equation}\label{hhh}
	P_{i+1} = U_{i+1}^{\top} H_{i+1} U_{i+1}.
	\end{equation}
	Inserting (\ref{equ:18}) into (\ref{hhh}), it is straightforward to show that the iteration (\ref{hhh}) on $P_i$ converges to the solution of (\ref{def:P}) as $i$ tends to infinity, i.e., $\lim_{i \rightarrow \infty}P_{i} = P$ with $P$ given by (\ref{def:P}). The iterations on $g_{i+1}$ and $z_{i+1}$ can be analogously derived. Thus, the parameters $\theta_{i}$ converge to a solution determined by (\ref{def:P}) and (\ref{def:g}). Specifically,
	\begin{equation}\notag
	\begin{aligned}
	&\lim\limits_{i \rightarrow \infty}H_{i} = \begin{bmatrix}
	\alpha A^{\top} P A+ Q & \alpha A^{\top} P B & \alpha A^{\top} P E \\
	&\alpha B^{\top} P B+ R & \alpha B^{\top} P E \\
	*&  & \alpha E^{\top} P E-\lambda I
	\end{bmatrix}\\
	&\lim\limits_{i \rightarrow \infty}G_{i}^{\top} = \left[\alpha g^{\top}A, \alpha g^{\top}B, \alpha g^{\top}E+2\lambda \bar{w}^{\top}  \right]\\
	&\lim\limits_{i \rightarrow \infty}s_{i} = \alpha z - \lambda \|\bar{w}\|^2,
	\end{aligned}
	\end{equation}
	where $P,g,z$ are given by (\ref{def:P}) and (\ref{def:g}). The proof is completed. \qed
\end{pf}

\section{Numerical examples}\label{sec:experiment}
This section validates the convergence of the proposed Q-learning method via a numerical example.

Consider a quadrotor that operates on a 2-D horizontal plane with the following dynamics
$$
	x_{k+1}=\begin{bmatrix}
		1 & 0 & T & 0 \\
		0 & 1 & 0 & T \\
		0 & 0 & 1 & 0 \\
		0 & 0 & 0 & 1
	\end{bmatrix} x_{k}+\begin{bmatrix}
		\frac{T^2}{2} & 0 \\
		0 & \frac{T^2}{2} \\
		T & 0 \\
		0 & T
	\end{bmatrix}\left(u_{k}+w_{k}\right),
$$
where $T=0.1$, $[x_{k,1},x_{k,2}]$ denotes the position coordinates, $[x_{k,3},x_{k,4}]$ is the corresponding velocities, and the input $u_k$ is the acceleration. The disturbance $w_{k,1}$ is subject to a Gaussian mixture distribution of $\mathcal{N}(1,0.1)$ and $\mathcal{N}(0.6,0.1)$ with equal weights. The independent noise $w_{k,2}$ is sampled from $\mathcal{N}(0,0.1)$. We set
$
Q = I~\text{and}~R = 0.2 I.
$
The discount factor $\alpha$ is set to $\alpha = 0.99$. The length of the sampled trajectory is set to $M = 900$. We randomly generate a set of samples $\{v^{j}\}_{j=1}^{10}$ from its groundtruth distribution and obtain the sample mean
$\bar{w} =[0.681, 0.132]^{\top}$. 

\begin{figure}
	\centerline{\includegraphics[width=70mm]{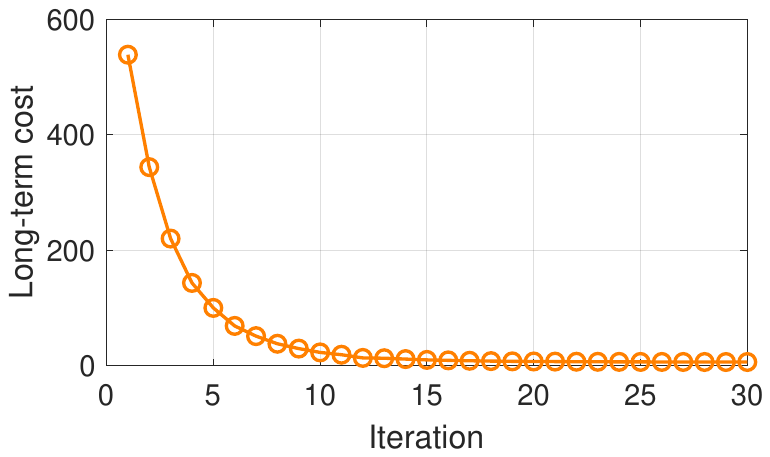}}
	\caption{The long-term cost $J_i$ converges within 30 iterations.}
	\label{pic:cost}
\end{figure}

By Assumption \ref{assum:lambda}, the penalty factor should satisfy $\lambda I -\alpha E^{\top} P E > 0$, which in our experiments is $\lambda > 0.22$. For demonstration, we select the penalty factor to $\lambda = 0.9$ and apply Algorithm \ref{alg:q_learning} to learn a minimax controller. We select the long-term cost at the $i$-th iteration
$$
J_i = \sum_{k=0}^{M-1} \alpha^{k} ((x_k^{i})^\top Qx_k^i + (u_k^{i})^\top Ru_k^i - \lambda \|w_k^i- \bar{w}\|^2 )
$$
as an indicator for convergence, where $u_k^{i} = K_ix_k^i+ r_i$ and $w_k^{i} = L_ix_k^i+ l_i$ are the pair of optimal solutions at the $i$-th iteration. The result is displayed in Fig. \ref{pic:cost}. As indicated by Theorem \ref{thm:convergence}, $J_i$ converges exponentially, which implies that the Algorithm \ref{alg:q_learning} returns an optimal Q-function.

\section{Conclusion}\label{sec:conclusion}
We have proposed a Q-learning algorithm with convergence guarantees to learn a minimax controller with robustness to disturbance distribution errors. Recent years have witnessed the tremendous success in the policy gradient methods, which can be used to solve the zero-sum game~\citep{zhang2019policy}. We shall study this topic in the future work.
\bibliographystyle{agsm}      
\bibliography{mybibfile}
\end{document}